\documentclass[sn-mathphys]{sn-jnl}
\usepackage{setspace}
\usepackage[utf8]{inputenc}
\usepackage[T1]{fontenc}
\usepackage{amsmath}
\usepackage{amsthm}
\usepackage{eucal}
\usepackage{amssymb}
\usepackage{mathrsfs}
\usepackage[short, c1]{optidef}
\usepackage{booktabs}
\usepackage{cuted}
\usepackage{textcomp}
\usepackage{tgpagella} 
\usepackage{natbib}
\usepackage{siunitx}
\usepackage{eurosym}
\usepackage{physics}
\usepackage{xcolor}
\usepackage{microtype}
\usepackage{graphicx}
\usepackage{caption,subcaption}

\usepackage{color}
\usepackage{ifthen}
\usepackage{soul}
    \definecolor{DarkGreen}{rgb}{0.2,0.5,0.2} 
\newcommand{\withColorMarking}{true} 
\ifthenelse{\equal{\withColorMarking}{true}}
  {  }
  {  }
\usepackage{array}
\usepackage{multirow}

\definecolor{DarkBlue}{rgb}{0.2,0.2,0.5} 

\hypersetup{pdftitle={},
            pdfauthor={},
            pdfkeywords={},
            pdfdisplaydoctitle,
            pdfstartview=Fit}

\begin{document}

\title{Federated Deep Learning in Electricity Forecasting: An MCDM Approach}
\author[1,2,3]{\fnm{Marco} \sur{Repetto}}\email{marco.repetto@skema.edu}
\author*[1]{\fnm{Davide} \sur{La Torre}}\email{davide.latorre@skema.edu}
\author[4]{\fnm{Muhammad} \sur{Tariq}}\email{m.tariq@adsm.ac.ae}

\affil*[1]{\orgname{SKEMA Business School, Université Côte d'Azur}, \city{Sophia-Antipolis}, \country{France}}
\affil[2]{\orgname{University of Milan-Bicocca}, \city{Milan}, \country{Italy}}
\affil[3]{\orgname{Siemens S.p.a.}, \city{Milan}, \country{Italy}}
\affil[4]{\orgname{Abu Dhabi School of Management}, \city{Abu Dhabi}, \country{UAE}}

\abstract{Large-scale data analysis is growing at an exponential rate as data proliferates in our societies. 
This abundance of data has the advantage of allowing the decision-maker to implement complex models in scenarios that were prohibitive before.
At the same time, such an amount of data requires a distributed thinking approach. In fact, Deep Learning models require plenty of resources, and distributed training is needed. This paper presents a multicriteria approach for distributed learning. Our approach uses the Weighted Goal Programming approach in its Chebyshev formulation to build an ensemble of decision rules that optimize aprioristically defined performance metrics. 
Such a formulation is beneficial because it is both model and metric agnostic and provides an interpretable output for the decision-maker.
We test our approach by showing a practical application in electricity demand forecasting. Our results suggest that when we allow for dataset split overlapping, the performances of our methodology are consistently above the baseline model trained on the whole dataset.}

\maketitle
\section{Introduction}

It is well know that Artificial Intelligence (AI) identifies in a broad sense the ability of a machine to learn from experience, to simulate the human intelligence, to adapt to new scenarios, and to get engaged in human-like activities. 
AI identifies an interdisciplinary area which includes computer science, robotics, engineering, mathematics. Over the years, it has made a rapid progress: it will contribute to the society transformation through the adoption of innovating technologies and creative intelligence and the large-scale implementation of AI in technologies such as IoT, smart speakers, chat-bots, cybersecurity, 3D printing, drones, face emotions analysis, sentiment analysis, natural language processing, and their applications to human resources, marketing, finance, and many others. 

With the term Machine learning (ML), instead, we identify a branch of AI in which algorithms are used to learn from data to make future decisions or predictions. ML algorithms are trained on past data in order to make future predictions or to support the decision making process.

Deep Learning (DL), instead, is a subset of ML and it includes a large family of ML methods and architectures based on Artificial Neural Networks (ANNs). It includes Deep Neural Networks, Deep Belief Networks, Deep Reinforcement Learning, Recurrent Neural Networks and Convolutional Neural Networks, to mention a few of them.
DL algorithms have been used in several applications including computer vision, speech recognition, natural language processing, bioinformatics, medical image analysis, and in most of these areas they have demonstrated to perform better than humans.
In the recent years DL has disrupted every application domain and it provides a robust, generalized, and scalable approach to work with different data types, including time-series data \cite{lecun2015deep, bengio2013representation, he2016deep, van2015lbann}. 

Fostered by the abundance of data, many recent DL architectures require a considerable amount of training.
At the same time, local regulations posed significant constraints in terms of data transmission in distributed systems \cite{ahmed2021}.
Because of this \cite{konecny2017} proposed the concept of Federated Learning (FL).  FL is a distributed learning methodology allowing model training on a large corpus of decentralized data \cite{bonawitz2019}. 

Recent applications of the FL focus primarily on the medical field, as particular attention has to be placed on patient information.
In \cite{jimenez-sanchez2021} the authors proposed an FL system based on memory-aware curriculum learning to solve the problem of class imbalance in the data of multiple institutions.
In \cite{beguier2021}, the authors proposed a sound FL methodology with Differential Privacy with an acceptable compromise between privacy and performance.
It is worth noting how both techniques may suffer from data leakage as by sharing gradients within nodes, some of the original information may be reconstructed \cite{zhu2020}.
Aside from the medical field, \cite{kholod2021} investigated FL applied on the Internet of Things (IoT) comparing different out-of-the-shelf open-source approaches.
In particular, the authors evaluated several aspects ranging from ease of use to performance in terms of implementation.

Aside from the performance, FL carries different concerns in terms of data leakage depending on the method utilized in the aggregation phase \cite{lim2021}.
Moreover, as pointed out by \cite{bagdasaryan2020} these systems can be subject to adversarial attacks.

Adversarial training is a fairly recent field in ML and DL. It is about the limitations of existing algorithms in their ability to correctly classify when they have been trained on adversarially datasets or to react against the effect of inputs purposely designed to fool a ML model.
A critical scenario in adversarial learning happens when training the same model from different datasets that show conflicting trends or behaviours. For instance, this is the case in financial forecasting when one tries to predict the future evolution of financial markets based on covid or pre-covid data.

Adversarial training impacts FL based on how the aggregation procedure is carried out.
As demonstrated in \cite{lim2021} in the case of FL in which gradients are shared, the attacker can leverage models gradients to either reconstruct datapoints or to inject an adversarial attack.
In particular, adversarial attacks may happen even in federated averaging as these aggregation techniques lack node accountability \cite{bagdasaryan2020}.

Motivated by the above considerations, in order to avoid the bias given by one dataset with respect to another, in this paper we propose an new approach based on Multiple Criteria Decision Making (MCDM).  The different loss functions obtained by training the same model on different datasets are combined in a unique framework by means of MCDM techniques. More in particular, we propose a Goal Programming (GP) model to deal with this complex scenario. Numerical experiments show the goodness of this approach and provide strategic managerial insights.

The intersections between MCDM and DL have not been explored yet and in the literature there are just a few contributions that try to combine MCDM and forecasting/machine learning techniques.
To mention a few of them, for instance in \cite{angelopoulos2019disaggregating} the authors analyze the relationship between times series and influential multiple criteria and they offer long-term electricity demand predictions in Greece. 
In \cite{stoilova2021novel}, instead, the authors propose a novel fuzzy multiple criteria time series modelling method based on fuzzy linear programming and sequential interactive techniques and they apply it to urban transportation planning. 

Our approach aligns with the FL paradigm as our dataset is heterogeneous and spread across many nodes with some computing capabilities.
Moreover, our approach is affected by the literature on adversarial training as models compete in the aggregation phase.
In our case, the number of models coincides with the number of nodes, but theoretically, this can be extended by allowing ensemble in each local node.
Our approach distinguishes from the previous ones as it is not model-dependent \cite{kholod2021}.
Moreover, our approach is different from model averaging as proposed by \cite{yager1988}, \cite{warnat-herresthal2021} and \cite{mcmahan2017} as we allow for nodes accountability and attribute nodes relevance based on their performance.

The paper is organized as follows. Section \ref{sec2} is devoted to the literature review on MCDM. In Section \ref{sec3} we present our main model. In Section \ref{sec4} we present numerical experiments and discuss the results of our approach in the case of electricity demand forecasting. Section \ref{sec5} concludes.

\section{Multiple Criteria Decision Making and Goal Programming} 
\label{sec2}

MCDM is a branch of Operations Research and Management Science. MCDM includes the set of methods and processes through which the concern for several conflicting criteria can be explicitly incorporated into the analytical process \cite{ehrgott2002}.
Several strategies, spanning from a priori to interactive ones, have been developed to address this issue.
Furthermore, these methods were widely used in a variety of fields, including economics, engineering, finance, and management \cite{colapinto2017multi}.
Historically MCDM's first roots are the ones laid by Pareto at the end of the 19th century. 
An MCDM problem is expressed mathematically as follows:

\begin{alignat}{3}
&\min_{\pmb{x}}     &\qquad& (f_1(\pmb{x}), f_2(\pmb{x}), \dots, f_p(\pmb{x}))\label{eq:optProb}\\
&\text{s.t.} &      & \pmb{x} \in S\label{eq:constraint1}
\end{alignat}  

where $f_i: {\rm I\!R}^n \rightarrow {\rm I\!R}$ is the ith objective and the vector $\pmb{x}$ contains the decionsion variabels that belong to the feasible set $S$.

Scalarization is a common approach for dealing with Multicriteria optimization problems. By scalarization, a vector optimization problem is converted into a single objective optimization problem by means of a weighted sum scalarization $ f(\pmb{x}) = \pmb{w}^\top \pmb{f}(\pmb{x}) = \sum_{i=1}^p w_i f_i(\pmb{x}) $ where $\pmb{w}$ is a vector weights. Each criterion is included in the scalarized version with different weights which express the relative importance of that criterion for the Decision Maker (DM). 
Another method that can be used to solve vector-valued problems is Goal Programming (or GP approach). The method was initially proposed by \cite{charnes1955optimal} and it aims at minimizing the deviations coming from the over or under achievement of different objectives.
The innovative idea behind this model is the determination of the aspiration levels of an objective function. This model does not try to find an optimal solution but an acceptable one (or satisfying one), and it tries to achieve the goals set by the DM rather than maximising or minimising the different criteria. Given a set of ideal goals $g_i$ , with $i = 1,..., p$, chosen by the DM, a weighted GP problem takes the form:

\begin{alignat}{3}
&\min_{\pmb{x},\delta^+_i,\delta^-_i}     &\qquad& \sum_{i=1}^{p} w^+_i \delta^+_i + w^-_i \delta^-_i \\
&\text{s.t.} &      & f_i(\pmb{x}) + \delta^-_i - \delta^+_i = g_i, \quad i = 1 \dots p \\
& &      & \pmb{x} \in S\label{eq:constraint2} & \\
& &      & \delta^-_i, \delta^+_i \geq 0, \quad i = 1 \dots p \label{eq:constraint3}
\end{alignat}  

where $\delta^+_i, \delta^-_i$ measures respectively the negative and the positive deviations to the optimal level $g_i$, whereas $w^+_i, w^-_i$ express the relative importance of the modeler to these goals.

\section{Model Formulation}
\label{sec3}

Before specifying our approach is worth defining some notation that will be used through the section and in the empirical application.
In our setup, we observe a tuple of the target variable and features over time, that is $\{(y_t, \boldsymbol{x}_t); t = 1 ... T\}$.
As our empirical application focuses on univariate time series forecasting, the features are the lagged values of the target variable. 
However, our approach is suitable for any use case in which a large dataset prevents any model estimation in a single machine or whether the dataset is distributed across many nodes not belonging to the same device.
Now let's assume that we slice the whole dataset into $K$ partitions such that $\sum_{k=1}^K S_k \geq T $.
In other words, the sum of the lengths of each partition is at least the length of the dataset itself. This formulation allows for overlapping in between partitions.
This is a recurrent situation in many distributed systems, as data redundancy is used to enforce the system's reliability.
In our empirical experiment, we considered the different degrees of overlapping ranging from no overlapping at all, that is, $\sum_{k=1}^K S_k  = T$ to full overlapping in which each partition includes the two adjacent ones. In other words $ \sum_{k=1}^K S_k \approx 3T$  

Given these premises, the proposed FL approach can be divided into three steps: local learning, federated validation-aggregation, and local updating.

In step 1, we train a model for each data partition, that is $f_k(w_k, \boldsymbol{x}_{k,t})$,  $\forall k \in K, \forall t \in S_k,$ where $S_k$ index the subset of each partition, we can define the predictions of such model as $\hat{y}_{k,t}$ where the prediction depends on both the model indexed with $k$ and the subset of features provided indexed with $t$.
In step 2, we solve the following Chebyshev Goal Programming \citep{flavell1976} problem to find an ensemble of models that jointly performs better than the sole models on all the dataset splits:

\begin{mini}<b>
{\alpha_j, \lambda}{\lambda}{}{}
\addConstraint{\delta_k^{+}}{\leq \lambda}{\forall k \in K}
\addConstraint{\sum_{j=1}^{K}\alpha_j L_j(y_{k,t}, \hat{y}_{k,t})- \delta_k^+ + \delta_k^-}{ =  L_k(y_{t}, \hat{y}_{k,t})}{\ \ \ \forall k \in K, \forall t \in S_k}
\addConstraint{\sum_{k=1}^K{\alpha_k} }{= 1}{}
\addConstraint{\alpha_k, \delta^{+}_k, \delta^{-}_k}{\geq 0}{\quad \forall k \in K}
\end{mini}

In step 3, with the optimal set of weights $\boldsymbol{\alpha}$ we predict the output as a convex combination of the $1...K$ trained models, that is: $$\hat{y_t} = \sum_{k=1}^K \alpha_k f_k(w_k, \boldsymbol{x}_{k,t}) $$.

The advantage of such an approach is that it is both model and metric agnostic.
More specifically, the approach is model agnostic in the sense that $f(\cdot)$ can be any predictive function, ranging from Deep ANNs to Gradient Boosted Machines. 
Moreover, the approach is metric agnostic as $L(\cdot)$ can be substituted with any, possibly non-differentiable, performance metric. An example in the case of classification may be to use accuracy or the Area Under the Receiving Operating Curve (AUROC). 
However, in such cases, as both the accuracy and AUROC need to be maximized, $\delta^-$ should be placed in the first soft constraint of the GP.

\section{Numerical Experiments and Discussion}
\label{sec4}

This section describes an application of our approach to the use case of electricity demand forecasting.
In particular, we employ a publicly available dataset, the 2017 Global Energy Forecasting Competition (GEFCom2017) dataset \citep{hong2019}. 
The dataset consists of information about the electricity load, holiday, and weather data.
Our application in electricity demand forecasting focuses only on electricity load, particularly on the total electricity load given by the sum of all the different New England zones.

\begin{figure}[t]
\begin{center}
\includegraphics[width=0.6\textwidth]{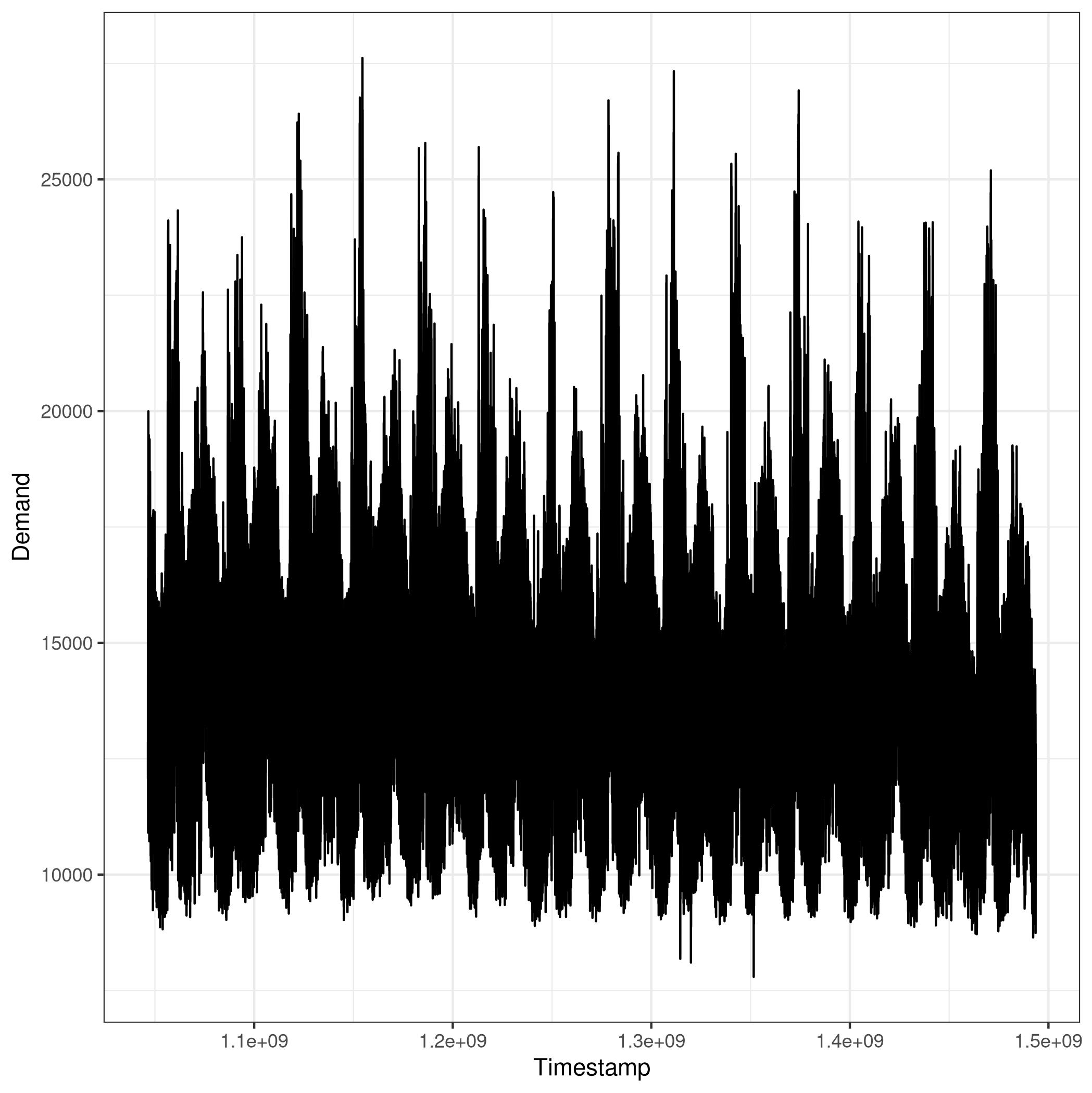}
\end{center}
\caption{Total electricity demand for the New England from 1st March 2003 to 30th April 2017 on a hourly basis.}
\label{total_ts}
\end{figure}

The plot in Figure \ref{total_ts} gives a perspective on the size of the time series that spans from 1st March 2003 to 30th April 2017 with a frequency of one hour.
In total, the time series consists of more than 100 thousand observations.

\subsection{Experimental setup}
The GEFCom2017 competition required only one month ahead forecast.
In our setup, we extend the forecasting horizon in order to compare our approach to the approach proposed by other authors, such as the one proposed in \cite{wang2020}.
Therefore we split the time series into training and test set.
The training set starts at the beginning of the dataset and ends on the 31st of December 2016.
The test set starts on the 1st of January 2017 and ends on the 30th of April 2017.
The forecasting horizon is then of 2879 steps.
To model this time series, we use the Feedforward Artificial Neural Network (FANN) in its most straightforward architecture, the Multilayer Perceptron (MLP).
The functional form of the MLP can be elicited as in \cite{arifovic2001}:
\begin{equation}
	\label{FANN_eqn}
	f(x) = \phi \left(\beta_0 + \sum_{j=1}^d\beta_jG \left(\gamma_{j0}+\sum_{i=1}^p \gamma_{ji}x_i 
	\right)\right)
\end{equation}
where G is the activation function, in our case $G(x) = \frac{1}{1+e^{-\alpha x}}$, $\beta$ and $\gamma$ represent weights and biases at each layer, whereas $\phi(\cdot)$ is the network output function that in our case is the identity function as our goal is to regress electricity demand and not time series classification.
To test our approach's robustness, we used different configurations in terms of split size and splits overlapping.
For the split size, we considered  10, 150, and 200, whereas, for the overlapping ratio, we considered complete, half, and no overlaps in between the splits.

\subsection{Evaluation metrics}
The evaluation metrics we used are: Root Mean Squared Error (RMSE), Mean Absolute Error (MAE), Mean Absolute Scaled Error (MASE) and Symmetric Mean Absolute Percentage Error (SMAPE).
RMSE is a standard metric used both in time series contexts and more broadly in any supervised learning setting. 
It is defined as:
 $$
 RMSE = \sqrt{\frac{1}{n} \sum_{i=1}^{n} (y_i - \hat{y}_i)^2}
 $$
the RMSE is a scale-dependent metric and, because it uses squared errors, is more sensitive to higher deviations than other measures such as MAE.
MAE, while being also a scale-dependent measure metric, is formulated as:
 $$
 MAE = \frac{1}{n} \sum_{i=1}^{n} \abs{y_i - \hat{y}_i}
 $$
Differently from the first two metrics, MASE is a scale-independent metric that is used chiefly in time series contexts, initially proposed by \cite{hyndman2006}, MASE is defined as follows:
 $$
  MASE \frac{\frac{1}{J}\sum_j\abs{y_j - \hat{y}_j}}
      {\displaystyle\frac{1}{T-1}\sum_{t=2}^T \abs{y_{t}-y_{t-1}}}.
$$
Firstly appeared in \cite{armstrong1985} SMAPE is also an accuracy metric used in time series. Contrary to the MASE, the SMAPE is a percentage based metric formulated as:
$$
SMAPE = \frac{1}{n} \sum_{i=1}^{n} \frac{\abs{y_i - \hat{y}_i}}{\frac{\abs{y_i} + \abs{\hat{y}_i}}{2}}
$$
Both SMAPE and MASE are used extensively in time series, and this is the reason behind their inclusion \citep{makridakis2020, shankar2020}.

\subsection{Approach results}
To test the performances of our approach, we treated both the number of splits and the overlapping ratio as hyperparameters.
A glimpse of the results is presented in Table \ref{tab}.

\begin{table}
\centering
\caption{Approach performances based on different parametrization of split size and overlapping ratio.}
\label{tab}
\resizebox{0.7\columnwidth}{!}{%
\begin{tabular}{@{}llllll@{}}
\toprule
\textbf{Split} & \textbf{Overlap} & \textbf{RMSE} & \textbf{MAE} & \textbf{MASE} & \textbf{SMAPE} \\ \midrule
10             & 1.0              & 1588.22      & 1256.06     & 2.66      & 0.10      \\
200            & 1.0              & 2979.64      & 2474.38     & 5.25      & 0.18      \\
150            & 1.0              & 3352.54      & 2640.00     & 5.60      & 0.19      \\
200            & 0.0              & 4083.87      & 3398.31     & 7.21      & 0.23      \\
10             & 0.5              & 4598.62      & 3548.89     & 7.53      & 0.24      \\
150            & 0.0              & 4545.16      & 3710.49     & 7.87      & 0.25      \\
200            & 0.5              & 4568.74      & 3747.89     & 7.95      & 0.25      \\
150            & 0.5              & 5105.29      & 3962.41     & 8.40      & 0.26      \\
10             & 0.0              & 9211.98      & 8755.59     & 18.57     & 0.49      \\ \bottomrule
\end{tabular}
}
\centering
\end{table}

\begin{figure}
\begin{center}
\includegraphics[width=0.6\textwidth]{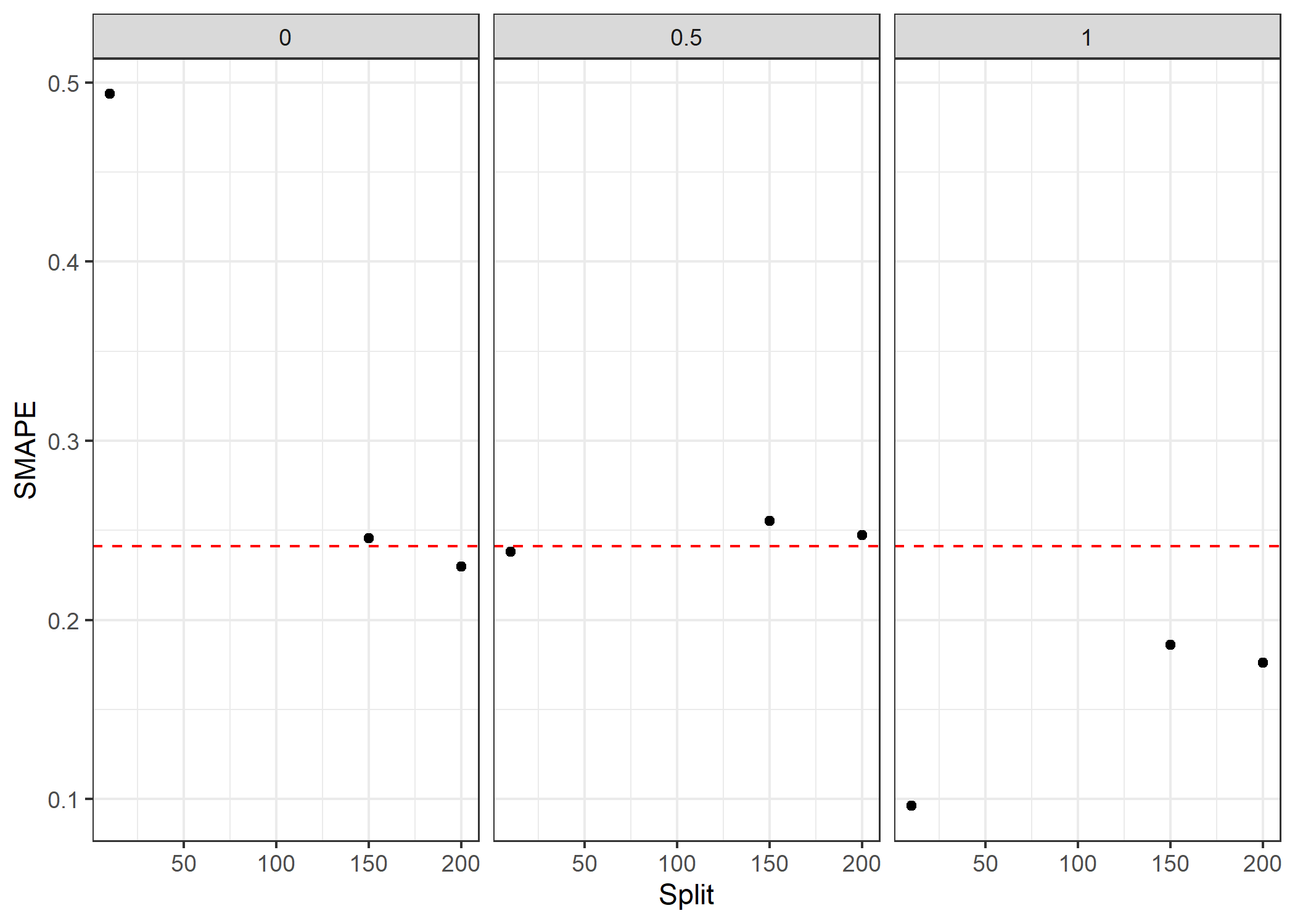}
\end{center}
\caption{Symmetric Mean Absolute Percentage Error on the test set under different split size and overlapping ratios.}
\label{robustnes}
\end{figure}

Figure \ref{robustnes} shows the different performances, measured with the SMAPE metric, of our approach and the sensitivity to hyperparameters change.
In particular, a low number of splits in combination with no overlaps results in a worse model than the baseline (dashed red line).
On the contrary, an increase in splits and overlapping results in an ensemble model on par with the benchmark.
It is worth noting that although the performance is similar between the two approaches in this setting, the advantage of our approach is that it allows for training in different machines and does not force the data to reside on the same device that is doing the model estimation.
Last, the more relevant result is that overlapping affects model performance positively. 
The result is far better than the benchmark, with a SMAPE of less than $0.1$ in the case of ten splits. 
More specifically, overlapping seems to offset the number of splits.
The potential of overlapping was already signaled by the literature \citep{scott1974analysis}.

\begin{figure}
\begin{center}
\includegraphics[width=0.6\textwidth]{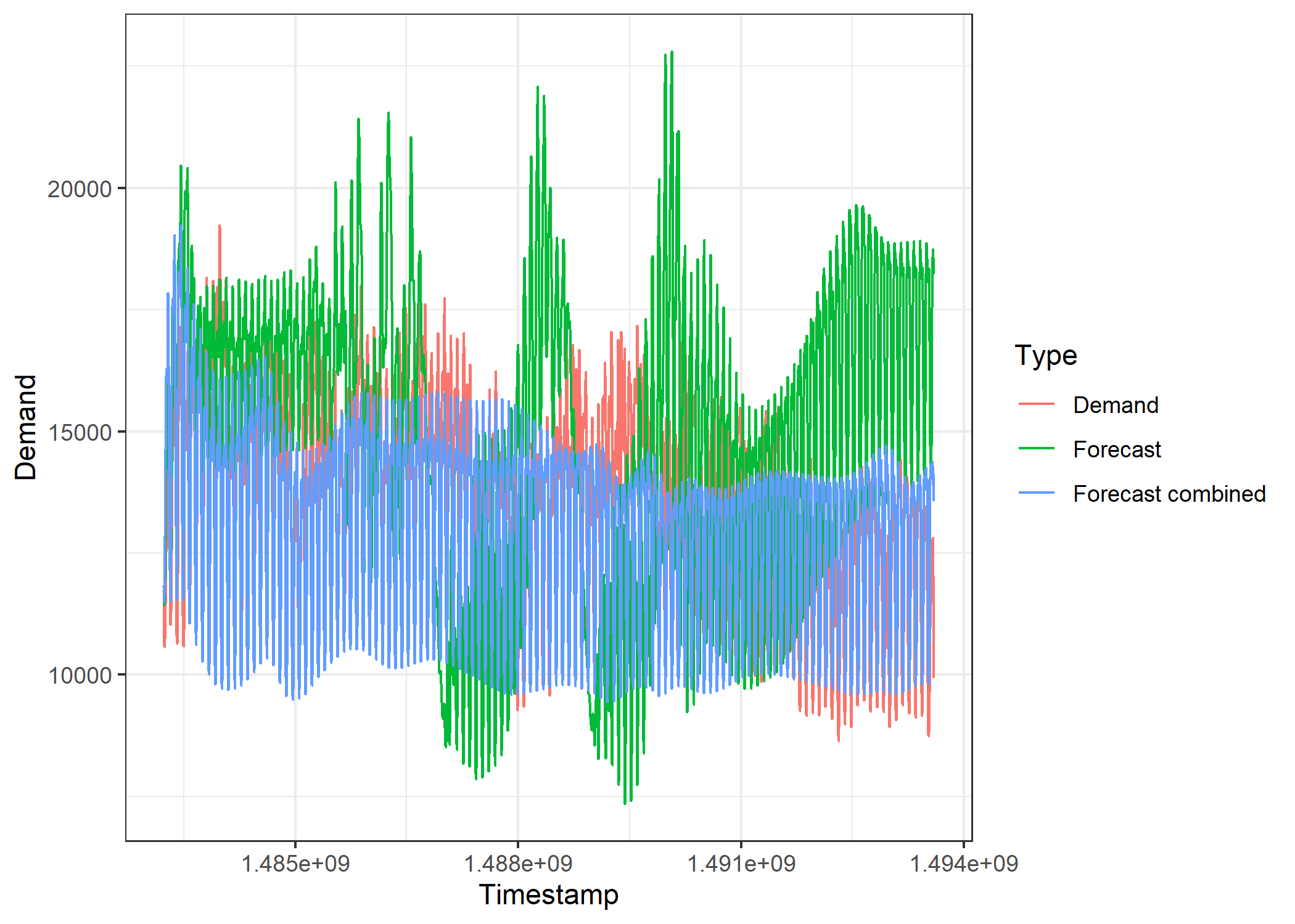}
\end{center}
\caption{Test set forecasts with respect to the true electricity demand (red line). The green line depicts the forecasts based on the model trained on the full dataset. The blue line depicts the forecasts based on the model trained with our approach with split size 10 and full overlapping.}
\label{forecast}
\end{figure}

This difference in performance is straightforward if we look at Figure \ref{forecast}.
In contrast to the model trained on the whole dataset, the combined forecast captures better the decreasing trend of the time series.
Moreover, the forecast is less volatile than the benchmark.

\begin{figure}
\begin{center}
\includegraphics[width=0.6\textwidth]{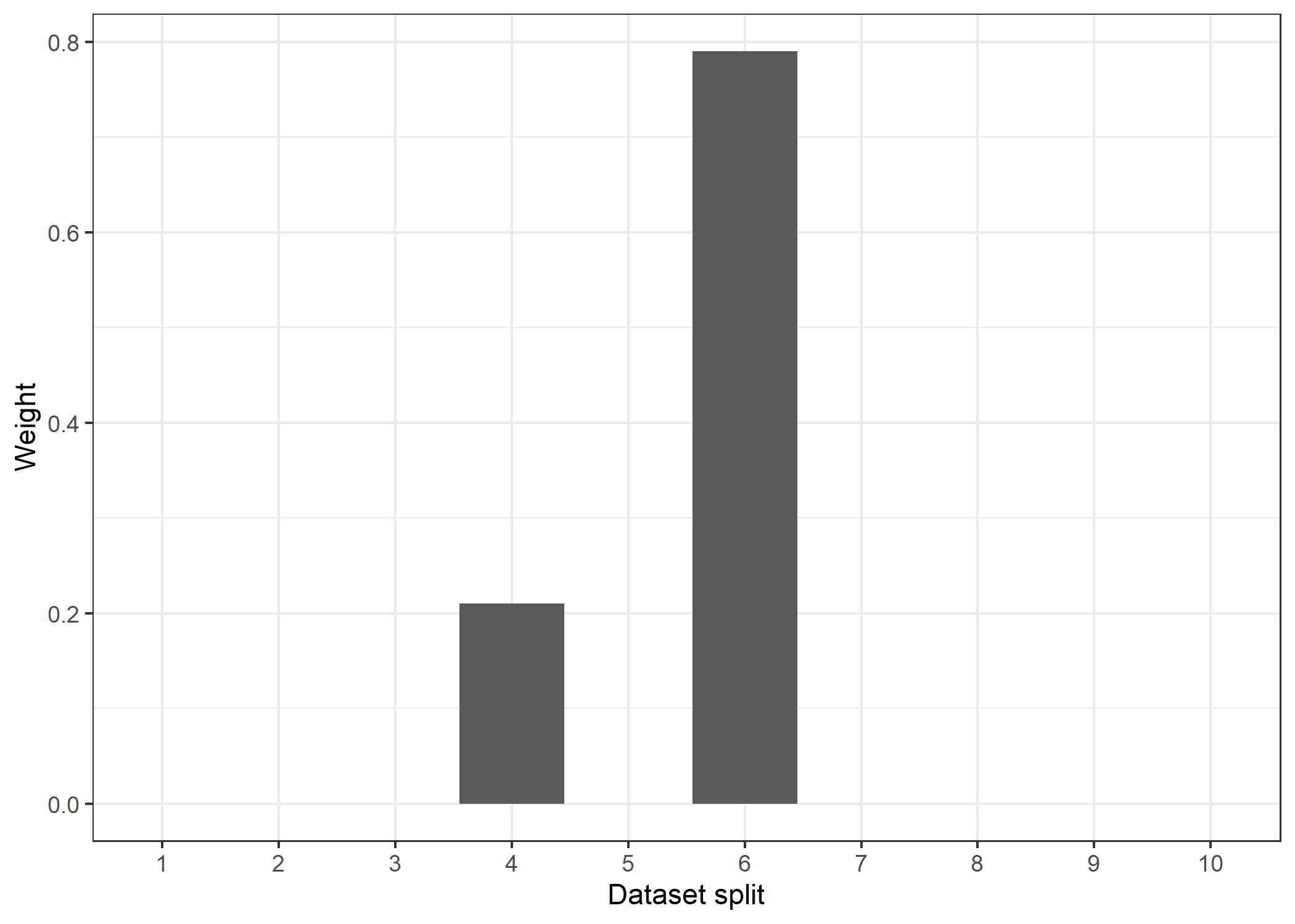}
\end{center}
\caption{Weights partitioning of the ensemble model based on split size 10 and complete overlapping.}
\label{weights}
\end{figure}

To investigate better the inner workings of our approach is also worth looking at the weights attributed by our multicriteria decision problem.
These weights are in Figure \ref{weights} and are attributed to the model trained on a particular data split.
Therefore, for example, the model trained on the first split, the most recent one, has zero weights and does not contribute to the final forecast.
On the contrary, the models trained on the fourth and sixth split are the only ones with non-zero weights.
This means that these models are the ones that score best in comparison with other models trained on different splits.

\begin{figure}
\begin{center}
\includegraphics[width=0.6\textwidth]{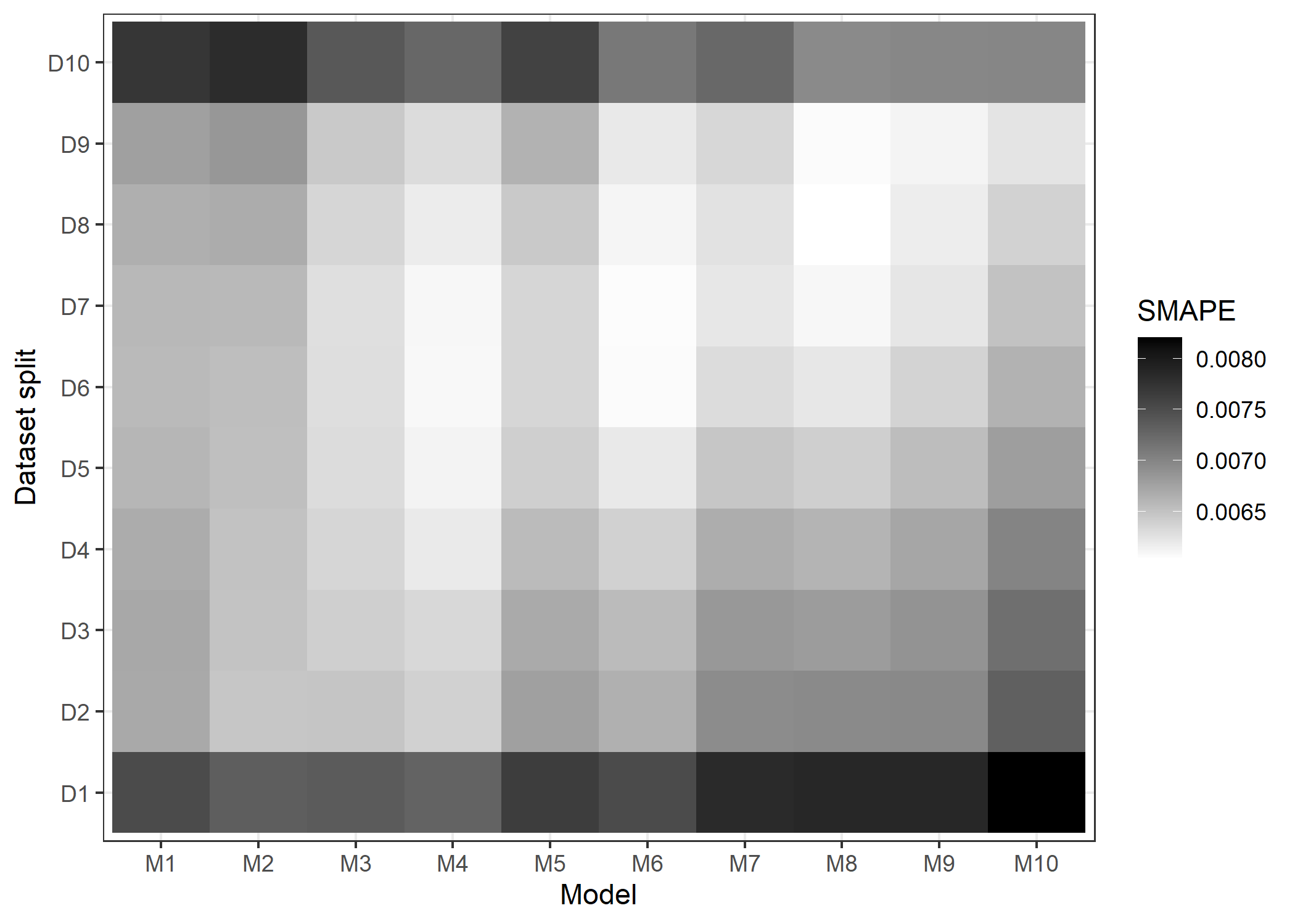}
\end{center}
\caption{Symmetric Mean Absolute Percentage Error on all the training data splits based on the model.}
\label{tile}
\end{figure}

These dynamics are particularly highlighted in Figure \ref{tile}.
This figure represents the SMAPE of each model applied to a different split.
The model that performs worst is the tenth model, and we can see that by looking at the last column, which is the one with the darkest cells.
Another relevant information is that every model performs poorly on both the first split and the last one.
What, instead, is clear is that both M4 and M6 perform better in almost all the splits.

\section{Conclusion}
\label{sec5}

The presence of adversarial training, biased datasets, and conflicting trends can generate problems during the machine learning process. This appears to be very crucial when thinking that splitting a dataset into training, test, and hold-out subsets is, most of the time, arbitrary and only led by the need to improve the algorithm performance and accuracy.  
On the top of this, dealing with huge amount of information can lead to long execution time of the gradient descent algorithm. 
In order to get efficient and prompt decisions it looks like of paramount importance to be able to split the original dataset into smaller chunks, training the same or different models on each of them separately, and finally recombine all of them in a unique framework. 
The advantages of this approach are twofold: from one side it allows to combine different models that are trained on different datasets without exchange of information during the gradient descent approach. This can be run on parallel architectures. From another the proposed GP formulation allows to recombine these models in a unique and efficient framework. Our numerical results confirm the goodness of this approach in the case of time-series analysis. Future avenues include the extensions of this framework by including different MCDM and GP approaches.

\bibliography{biblio}

\end{document}